\documentclass[10pt]{amsart}
\usepackage{amssymb}
\usepackage{amsmath,amssymb,amsfonts,amsthm,graphics,
latexsym, amscd, amsfonts, epsfig, eepic,epic}
\usepackage{mathrsfs}
\usepackage{color}
\usepackage{eucal}%    caligraphic-euler fonts: \mathcal{ }
\usepackage{eufrak}%   frak-euler        fonts: \mathfrak{ }
\usepackage[all]{xypic}
\usepackage{xspace}

\theoremstyle{plain}
\newtheorem{theorem}{Theorem}

\newtheorem{lemma}{Lemma}

\theoremstyle{definition}

\theoremstyle{example}

\theoremstyle{remark}

\numberwithin{equation}{section}

\begin{document}
%%%
%%%
%%%%%%%%%%%%%%%%%%%%%%%%%%%%%%%%%%%%%%%%%%%%%%%%%%%%%%%%%%%%%%%%%%%%%%%%%%
%%
%%%
\title[A combinatorial framework for RNA tertiary interaction]
      {A combinatorial framework for RNA tertiary interaction}
\author{Jing Qin and Christian M. Reidys$^{\,\star}$}
\address{Center for Combinatorics, LPMC-TJKLC %XXX%
           \\
         Nankai University  \\
         Tianjin 300071\\
         P.R.~China\\
         Phone: *86-22-2350-6800\\
         Fax:   *86-22-2350-9272}
\email{reidys@nankai.edu.cn}%XXXX
\thanks{}
\keywords{RNA structure, pseudoknot, base triple, asymptotic
enumeration, kernel method, partition, crossing, tangled-diagram}
\date{October, 2007}
\begin{abstract}
In this paper we show how to express RNA tertiary interactions via
the concepts of tangled diagrams. Tangled diagrams allow to
formulate RNA base triples and pseudoknot-interactions and to
control the maximum number of mutually crossing arcs. In particular
we study two subsets of tangled diagrams: $3$-noncrossing
tangled-diagrams with $\ell$ vertices of degree two and $2$-regular,
$3$-noncrossing partitions (i.e.~without arcs of the form
$(i,i+1)$). Our main result is an asymptotic formula for the number
of $2$-regular, $3$-noncrossing partitions, denoted by $p_{3,2}(n)$,
$3$-noncrossing partitions over $[n]$. The asymptotic formula is
derived by the analytic theory of singular difference equations due
to Birkhoff-Trjitzinsky. Explicitly, we prove the formula
$p_{3,2}(n+1)\sim K \ 8^{n}n^{-7}(1+c_{1}/n+c_{2}/n^2+c_3/n^3)$
where $K,c_i$, $i=1,2,3$ are constants.
\end{abstract}
\maketitle
{{\small
%\tableofcontents
}}

%%%
%%%%%%%%%%%%%%%%%%%%%%%%%%%%%%%%%%%%%%%%%%%%%%%%%%%%%%%%%%%%%%%%%%%%%%%%
%%%

\section{Introduction}\label{S:intro}

%%%
%%%%%%%%%%%%%%%%%%%%%%%%%%%%%%%%%%%%%%%%%%%%%%%%%%%%%%%%%%%%%%%%%%%%%%%%
%%%

It is well-known that the functional repertoire of RNA is closely
related to the variety of its shapes. Therefore it is of utmost
importance to understand the structural ``language'' of RNA as this
will eventually allow for fast folding, identification and discovery
of new RNA functionalities. Studies of RNA structural motifs at high
resolution by NMR and X-ray crystallographic methods provided
insight into the fundamental forces that give rise to the unique
structural characteristics of RNA. Non-Watson-Crick
purine-pyrimidine, purine-purine, and pyrimidine-pyrimidine base
pairing, as well as base-phosphate and base-ribose hydrogen bonding,
are known to be important forces for folding and stabilizing RNA
structures \cite{Shen}. For RNA pseudoknots (viewed as interactions
between unpaired bases) combinatorial abstractions have led to new
interpretations, generating functions and enumeration results.
Although far from having a complete understanding of RNA pseudoknots
conceptual progress has been made in identifying the right concepts
of, for instance crossing-complexity, which have direct implications
for novel RNA pseudoknot folding algorithms. In this paper we build
on the concepts derived in the context of RNA pseudoknots.

Before we begin by giving some background on RNA structure, let us
remark why ``combinatorial frameworks'' are of central importance
for any prediction algorithm. The above mentioned language of RNA is
tantamount to uniquely specifying each element of the variety of
shapes. Any prediction involves at some point a search through
configurations and has to make sure that shapes are, for instance,
not counted multiple times. The enumeration of the combinatorial
class and analysis of its mathematical structure are of fundamental
importance for designing such a search procedure. The primary
sequence of an RNA molecule is its sequence of nucleotides {\bf A},
{\bf G}, {\bf U} and {\bf C} together with the Watson-Crick ({\bf
A-U}, {\bf U-A}, {\bf G-C},{\bf C-G}) and ({\bf U-G}, {\bf G-U})
base pairings. Single stranded RNA molecules form helical structures
whose bonds satisfy the above base pairing rules and which, in many
cases, determine their function. Due to the biochemistry of the base
pairs stacked base pairs, i.e.~arcs of the form $(i,j),(i-1,j+1)$
have typically a lower minimum free energy than crossing arcs. Base
stacking is as important in determining RNA conformations as
hydrogen bonding interactions. With the noncanonical interactions,
many single-stranded loop regions such as hairpin loops, bulge
loops, and internal loops fold into well-defined secondary
structures. The prediction of RNA secondary structure is of
complexity $O(n^3)$ in time and $O(n^2)$ in space for a sequence of
length $n$ \cite{Zuker:81, Zuker:84} which is result from the fact
that no two bonds can cross.
%%%
%%%%%%%%%%%%%%%%%%%%%%%% Figures ex1 %%%%%%%%%%%%%%%%%%%%%%%%%%%%%%%%%%%%%%
%%%
\begin{figure}[ht]
\centerline{%
\epsfig{file=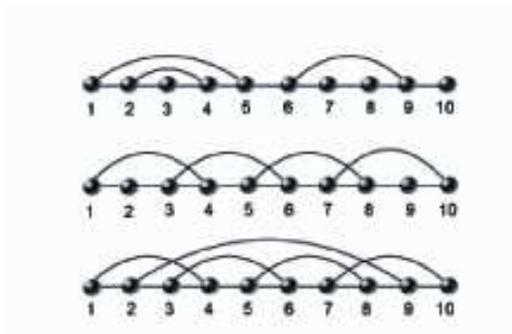,width=0.5\textwidth}\hskip15pt
 }
\caption{\small The idea behind the notion of $3$-noncrossing RNA structures.
(a) secondary structure (with isolated labels $3,7,8,10$),
(b) bi-secondary structure \cite{Stadler:99}, $2,9$ being isolated
(c) $3$-noncrossing structure, which is
    {\it not} a bi-secondary structure. In fact, this is {\it the} smallest
    $3$-noncrossing RNA structure which is not a bi-secondary structure.
}
\label{F:1}
\end{figure}
%%%
%%%%%%%%%%%%%%%%%%%%%%%%%%%%%%%%%%%%%%%%%%%%%%%%%%%%%%%%%%%%%%%%%%%%%%%%
%%%

While the concept of secondary structure is of fundamental
importance, it is well-known that there exist additional types of
nucleotide interactions \cite{Science:05}. These bonds are called
pseudoknots \cite{Westhof} and occur in functional RNA (RNAseP
\cite{Loria}), ribosomal RNA \cite{Konings} and are conserved in the
catalytic core of group I introns. Stadler {\it et al.}
\cite{Stadler:99} suggested a class of RNA pseudoknots called
bi-secondary structures which are essentially ``superpositions'' of
the arcs of two ``secondary structures'' and accordingly generalize
from outer-planar to planar graphs, see Figure~\ref{F:1}. Prediction
algorithms for RNA pseudoknot structures are much harder to derive
since there exists no {\it a priori} recursion and the subadditivity
of local solutions is not guaranteed. The key for enumerating RNA
pseudoknot structures is their categorization in terms of the
maximal size of sets of mutually crossing bonds
\cite{Reidys:07pseu}, i.e.~the notion of $k$-noncrossing structures.
To be precise, it is the inherent {\it locality} of the property
``$k$-noncrossing'' that allows for their enumeration by lattice
paths. The diagram representation  of a structure illustrates what
$k$-noncrossing means, see Figure~\ref{F:1}. In a diagram all
nucleotides are drawn horizontally and the backbone bonds are
ignored, then all bonds are drawn as arcs in the upper half-plane.
The number of $3$-noncrossing RNA structures satisfies ${\sf S}_3(n)
\sim \frac{10.4724\cdot
4!}{n(n-1)\dots(n-4)}\left(\frac{5+\sqrt{21}}{2}\right)^n$
\cite{Reidys:07asym}, however, it is not the exponential growth rate
of ($\frac{5+ \sqrt{21}}{2}$) but the inherent non-recursiveness
which makes the prediction difficult.

%%%
%%%%%%%%%%%%%%%%%%%%%%%%%%%%%%%%%%%%%%%%%%%%%%%%%%%%%%%%%%%%%%%%%%%%%%%%%%
%%%
\begin{figure}[ht]
\centerline{%
\epsfig{file=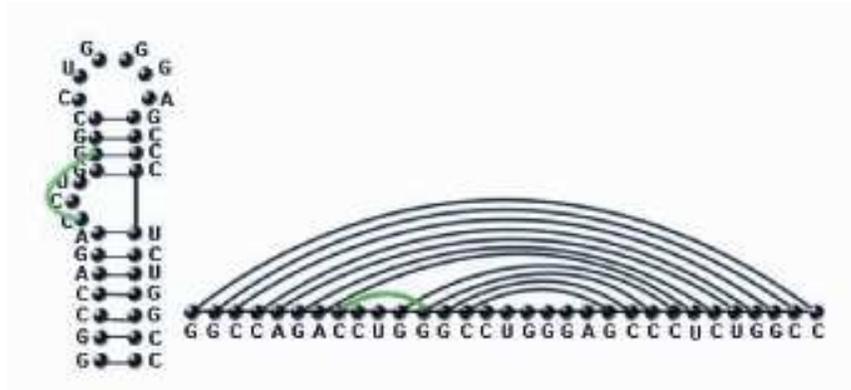,width=0.8\textwidth}\hskip15pt
 }
\caption{\small HIV-2 TAR, \cite{NMR}. In HIV-2 TAR we have a
$(\textbf{C}38$-$\textbf{G}27)\cdot\textbf{C}23^{+}$ triple mutant.
Improved NMR spectral properties of HIV-2 TAR allowed the
observation of the \textbf{C}23 amino and imino protons, providing
direct evidence of hydrogen bonding interaction. The tertiary
interaction is a tangled-diagram of with one vertex of degree two.}
\label{F:2}
\end{figure}
%%%
%%%
%%%%%%%%%%%%%%%%%%%%%%%%%%%%%%%%%%%%%%%%%%%%%%%%%%%%%%%%%%%%%%%%%%%%%%%%%%
%%%

\begin{figure}[ht]
\centerline{%
\epsfig{file=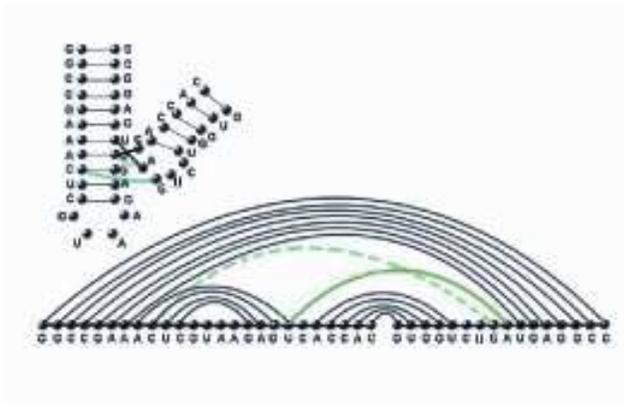,width=0.6\textwidth}\hskip15pt
 }
\caption{\small Diagram representation of the hammerhead ribozyme
\cite{Batey}, which can be represented as a tangled-diagrams  with
two vertices of degree two. The gap after $\textbf{C}$25 indicates
that some nucleotides are omitted, which are involved in an
unrelated structural motif.
 }
\label{F:3}
\end{figure}

A first step towards RNA-tertiary structures beyond pseudoknot
interactions consists in considering single strands interacting with
helical regions by forming tertiary contacts with base-paired
nucleotides of the helices. Nucleotide triples occur when
single-stranded nucleotides form hydrogen bonds with nucleotides
that are already base paired. This hydrogen bonds can involve bases,
sugars and phosphates. These interactions function to orient regions
of secondary structures in large RNA molecules and to stabilize RNA
three-dimensional structures. Base triples are a special case of
nucleotide triple interactions in which base-base hydrogen bonding
occurs. Single-stranded nucleotides can interact with base paired
nucleotides via either the major groove or the minor groove of
duplex regions. Nucleotide triples have been shown or proposed to
form at junctions of coaxially stacked RNA helices that have
adjacent single-stranded regions \cite{Shen,Batey}. Several major
groove triples are present in tRNA where they function to stabilize
its L-shaped three-dimensional structure. These interactions require
to consider tangled diagrams \cite{Reidys:07vac}, i.e.~diagrams with
vertices of degree $\le 2$ which exhibit a variety of arc
configurations, see Section~\ref{S:vac}. This variety is motivated
from nucleotide interactions observed in RNA structures. In
Figure~\ref{F:2} we show the HIV-2 TAR
$(\textbf{C}38$-$\textbf{G}27)\cdot\textbf{C}23^{+}$ triple mutant
structure as a tangled-diagram. Let us next have a closer look at
the hammerhead structure-motif \cite{Batey} in Figure~\ref{F:3}.
Comparing Figure~\ref{F:2} with Figure~\ref{F:3} reveals one feature
of the hammerhead motif. It exhibits a lefthand-endpoint of degree
$2$ (incident to the dashed arc) while all other vertices of degree
$2$ are left-and righthand-endpoints. These two examples indicate
that the majority of the bonds is organized in helical regions,
where Watson-Crick and {\bf G-U}({\bf U-G}) base pairs are stacked,
additional stacks can be realized forming pseudoknots.

%%%
%%%%%%%%%%%%%%%%%%%%%%%%%%%%%%%%%%%%%%%%%%%%%%%%%%%%%%%%%%%%%%%%%%%%%%%%%%
%%%
\begin{figure}[ht]
\centerline{%
\epsfig{file=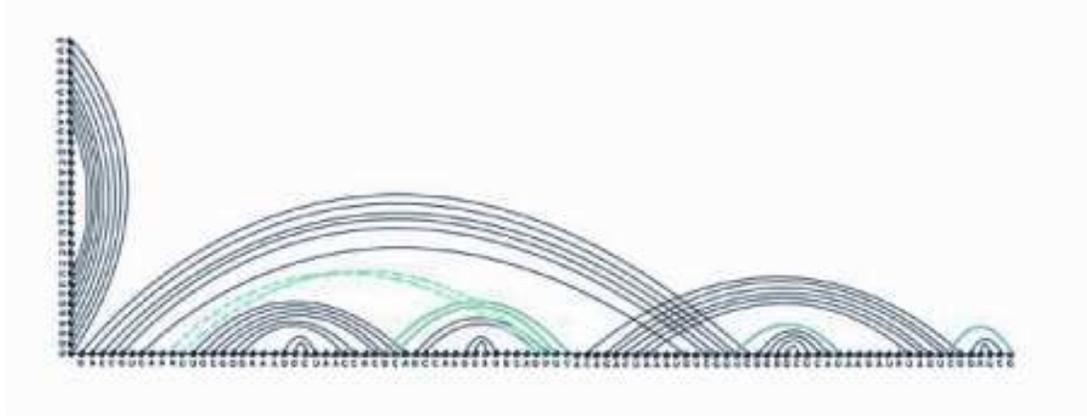,width=1.0\textwidth}\hskip15pt
 }
\caption{\small Catalytic core region of the group I self-splicing
intron \cite{Chastain} corresponds to a tangled-diagram with six
vertices of degree two. The gaps after $\textbf{G}54, \textbf{U}72,
\textbf{G}103$ and $\textbf{A}112$ indicate that some nucleotides
are omitted which are involved in an unrelated structural motif.
 }
\label{F:4}
\end{figure}

%%%
Finally in Figure~\ref{F:4} we display the catalytic core region of
the group I self-splicing intron \cite{Chastain}. In order to
express tertiary interactions we consider tangled-diagrams
introduced in \cite{Reidys:07vac}, which capture the nucleotide
interactions relevant for the tertiary structure of the molecule
\cite{Batey}.
%%%
%%%%%%%%%%%%%%%%%%%%%%%%%%%%%%%%%%%%%%%%%%%%%%%%%%%%%%%%%%%%%%%%%%%%%%%%%%
%%%

We will discuss two combinatorial frameworks arising from
tangled-diagrams \cite{Reidys:07vac}, both being suited for
expressing RNA tertiary interactions. The first is the set of
tangled-diagrams with fixed number of vertices of degree $2$ and the
second the set of $2$-regular $k$-noncrossing partitions. While the
former can easily be enumerated the latter requires more work.
$2$-regular $k$-noncrossing partitions evade lattice path
enumeration due to their inherent asymmetry (lacking arcs of length
$1$). The ``straightforward'' ansatz via Inclusion-exclusion applied
to the set of {\it all} $k$-noncrossing partitions revealed a
connection between seemingly unrelated combinatorial objects:
partitions and enhanced partitions, enumerated by Bousquet-M\'elou
and Xin \cite{MIRXIN,Reidys:07dual}. In Lemma~\ref{L:bijection}
\cite{Reidys:07dual} we show how this relation can be used to obtain
the enumeration. Subsequently, we prove the following a simple
formula for the numbers of $2$-regular $k$-noncrossing partitions
\begin{equation}
p_{3,2}(n+1)\sim K \ 8^{n}n^{-7}(1+c_{1}/n+c_{2}/n^2+c_3/n^3) \ ,
\end{equation}
where $K=6686.408973$, $c_1=-28,\ c_2=455.77778$ and
$c_3=-5651.160494.$ As for the quality of approximation we present
the sub-exponential factors in the table below, where
$g(n)=Kn^{-7}(1+c_1/n+c_2/n^2+c_3/n^3)$.
\begin{center}
\begin{small}
\begin{tabular}{|c|c|c|c|c|c|}
\hline
\multicolumn{6}{c}{\textbf{The Sub-exponential Factor}}\\
  \hline
  $n$ & $\rho_{3}(n)/8^n$ & $g(n)$ & $n$ & $\rho_{3}(n)/8^n$ & $g(n)$ \\
  21 & $1.479\times 10^{-6}$ & $1.726\times 10^{-7}$  & 81 & $2.270\times 10^{-10}$ & $2.264\times 10^{-10}$  \\\
  31 & $1.283\times 10^{-7}$ & $1.112\times 10^{-7}$  &91 & $1.033\times10^{-10}$ & $1.031\times 10^{-10}$ \\
  41 & $2.104\times10^{-8}$ & $2.026\times 10^{-8}$ &101 & $5.088\times10^{-11}$ & $5.081\times 10^{-11}$ \\
  51 & $5.011\times10^{-9}$ & $4.939\times 10^{-9}$  &  501 & $8.100\times10^{-16}$ & $8.095\times 10^{-15}$  \\
  61 & $1.524\times10^{-9}$ & $1.514\times 10^{-9}$  &1001 & $6.507\times 10^{-18}$ & $6.502\times 10^{-18}$  \\
  71 & $5.514\times10^{-10}$ & $5.493\times 10^{-10}$ &10001 & $6.672\times 10^{-25}$ & $6.668\times 10^{-25}$  \\
  \hline
\end{tabular}
\end{small}
\end{center}
Our analysis is based on the theory of Birkhoff-Trjitzinsky,
which seems to be somewhat overlooked. While the two original papers
\cite{Birkhoff,T:wimp} are hard to read, the paper of Wimp and
Zeilberger \cite{wimp} provides a good introduction and shows via
various examples of how to apply the theory. Since the method (if it
applies) is quite powerful we give an overview of the analytic
theory of singular difference equations in the Appendix.
%%%
%%%%%%%%%%%%%%%%%%%%%%%%%%%%%%%%%%%%%%%%%%%%%%%%%%%%%%%%%%%%%%%%%%%%%%%%%
%%%
\section{Vacillating tableaux and tangled-diagrams}\label{S:vac}
%%%
%%%%%%%%%%%%%%%%%%%%%%%%%%%%%%%%%%%%%%%%%%%%%%%%%%%%%%%%%%%%%%%%%%%%%%%%%
\subsection{Tangled-diagrams}
\begin{figure}[ht]
\scalebox{0.6}[0.6]{\includegraphics*[30,630][610,830]{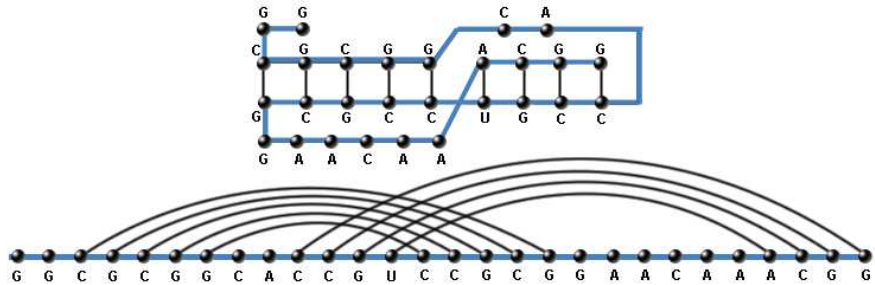}}
\caption{\small{We denote the backbone by the blue line and bonds by
black lines.}}
\label{F:pseudo}
\end{figure}

\begin{figure}[ht]
\scalebox{0.6}[0.6]{\includegraphics*[20,740][600,810]{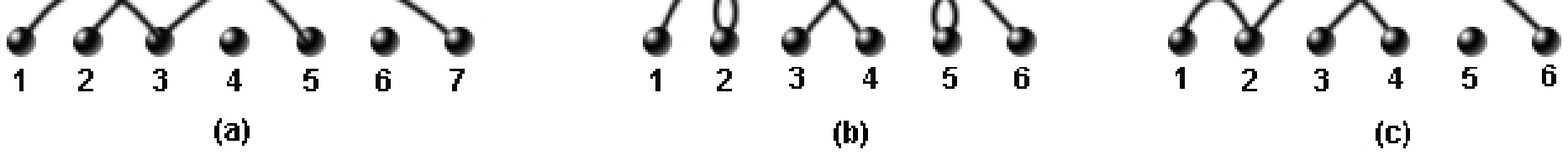}}
\caption{}\label{F:pli}
\end{figure}

A tangled-diagram over $[n]$ is a triple of sets $(V,E,F)$, where
$V$ is a finite non-empty set of $n$ elements called vertices, $E$
is a set of unordered pairs of vertices called arcs and $F$ is the
flag set whose elements are the 2-degree points such that they are
the ends of two crossing arcs, represented by drawing its vertices
in a horizontal line and its arcs $(i,j)$ in the upper halfplane
with the following basic configurations and the isolated points
%%%%%%%%%%%%%%%%%%%%%%%%%%%%%%%%%%%%%%%%%%%%%%%%%%%%%%%%%%%%%%%%%%5
%%%
\begin{center}\label{F:pb}
\scalebox{0.4}[0.4]{\includegraphics*[30,650][600,840]{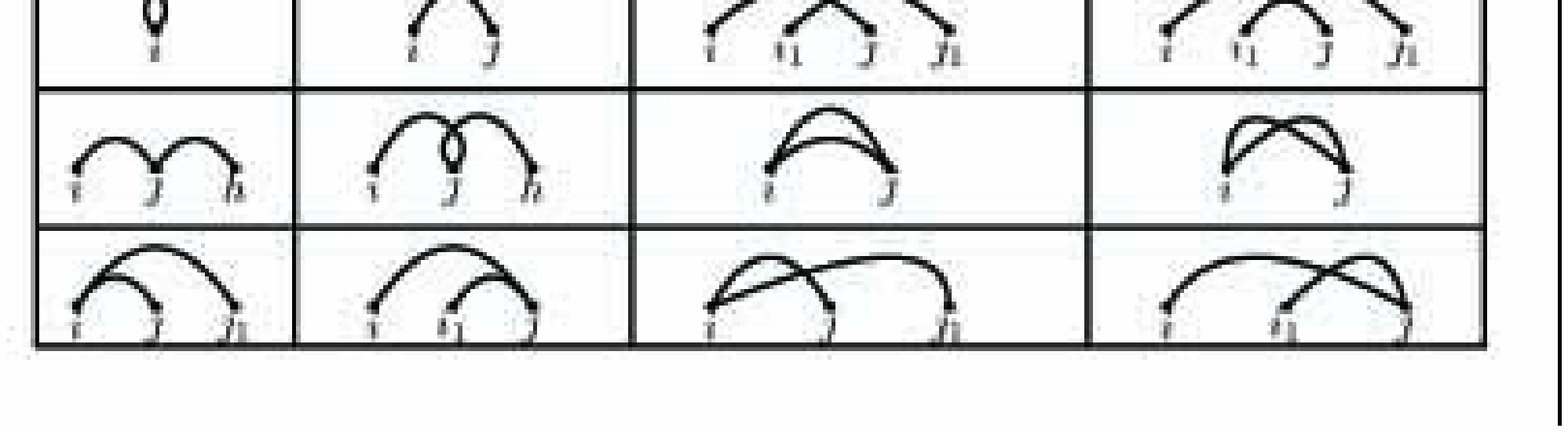}}
\end{center}
%%%
%%%%%%%%%%%%%%%%%%%%%%%%%%%%%%%%%%%%%%%%%%%%%%%%%%%%%%%%%%%%%%%%%%
%%%
Composing these motifs we obtain a tangled-diagram, for instance,
the tangled-diagram
%%%
%%%%%%%%%%%%%%%%%%%%%%%%%%%%%%%%%%%%%%%%%%%%%%%%%%%%%%%%%%%%%%%%%%5
%%%
\begin{center}\label{F:pb}
\scalebox{0.75}[0.75]{\includegraphics*[20,710][600,780]{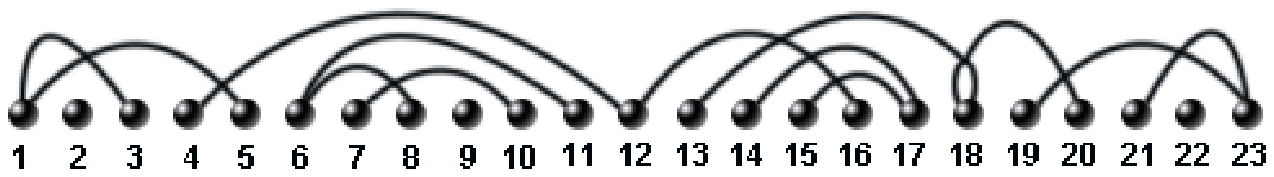}}
\end{center}
%%%
%%%%%%%%%%%%%%%%%%%%%%%%%%%%%%%%%%%%%%%%%%%%%%%%%%%%%%%%%%%%%%%%%%
%%%
has $V=[23]$ and $F=\{1,18,23\}$. Let us introduce several important
subclasses of $3$-noncrossing tangled-diagrams:\\
{\sf(1) $3$-noncrossing matchings with isolated points} are
$3$-noncrossing tangled-diagram in which each vertex has degree at
most $1$. For instance, RNA pseudoknot structures
\cite{Reidys:07pseu} are $3$-noncrossing matchings with isolated
points, see Figure \ref{F:pseudo}.
{\sf(2) $2$-regular, $3$-noncrossing partitions.} A partition
corresponds to a tangled-diagram in which any vertex of degree two,
$j$, is incident to the arcs $(i,j)$ and $(j,s)$, where $i<j<s$, for
instance, see Figure \ref{F:4} and Figure \ref{F:pli}, (a).
Partitions without arcs of the form $(i,i+1)$ are called
$2$-regular, partitions. {\sf(3) $3$-noncrossing braids without
isolated points} are tangled-diagrams in which all vertices, $j$ of
degree two are either incident to loops $(j,j)$ or crossing arcs
$(i,j)$ and $(j,h)$, where $i<j<h$, see Figure \ref{F:pli}, (b).
%%%
%%%%%%%%%%%%%%%%%%%%%%%%%%%%%%%%%%%%%%%%%%%%%%%%%%%%%%%%%%%%%%%%%%%%
%%%

%%%
%%%%%%%%%%%%%%%%%%%%%%%%%%%%%%%%%%%%%%%%%%%%%%%%%%%%%%%%%%%%%%%%%%%%
%%%
{\sf(4) $3$-noncrossing diagrams with $\ell$ vertices of degree 2.}
Figure~\ref{F:2}, Figure~\ref{F:3} and Figure~\ref{F:4} are
$3$-noncrossing tangled-diagrams with $\ell=1,2,6$ vertices of
degree $2$. The following tangle-diagram shows all $4$ basic types
of degree $2$ vertices in tangled diagrams.
%%%
%%%%%%%%%%%%%%%%%%%%%%%%%%%%%%%%%%%%%%%%%%%%%%%%%%%%%%%%%%%%%%%%%%5
%%%
\begin{center}
\scalebox{0.3}[0.3]{\includegraphics*[20,670][600,840]{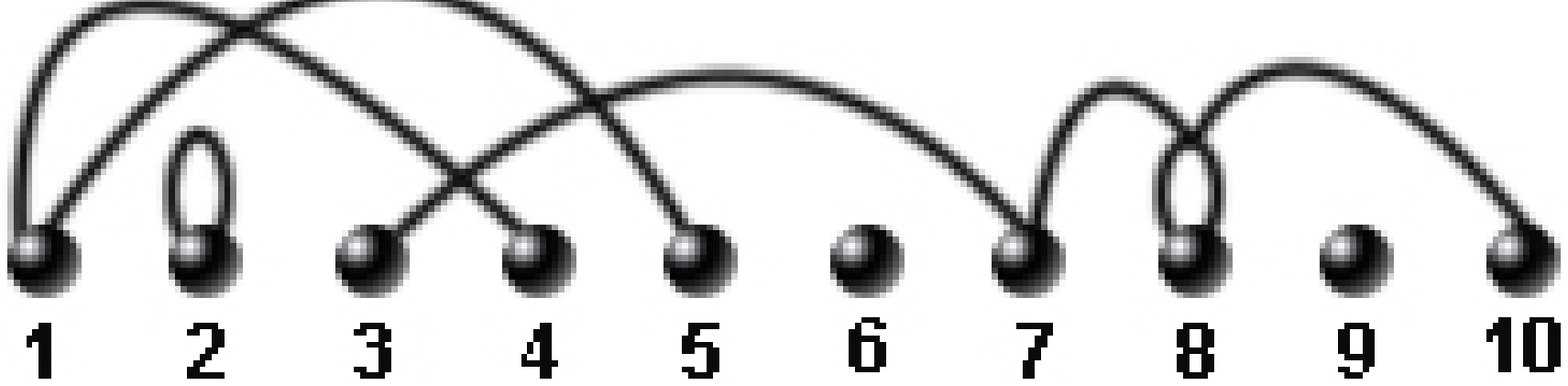}}
\label{F:diagram}
\end{center}
%%%
%%%%%%%%%%%%%%%%%%%%%%%%%%%%%%%%%%%%%%%%%%%%%%%%%%%%%%%%%%%%%%%%%%
%%%
In the following, we study the subclasses {\sf(2)} and {\sf(4)}
since they represent a natural framework for RNA tertiary
interactions. It turns out that {\sf (3)} is of importance since it
facilitates the enumeration of {\sf(2)}. To be precise it is shown
in \cite{Reidys:07dual} that there is a duality between
$k$-noncrossing braids without isolated points and $2$-regular
$k$-noncrossing partitions.

Having introduced the combinatorial framework, one key question is how to
enumerate the subclasses {\sf (2)} and {\sf (4)}. The enumeration is
facilitated via a bijection between the tangled-diagrams and certain
lattice paths. To derive the latter a bijection between tangled-diagrams
and (generalized) vacillating tableaux is constructed. It is then easy to
see that vacillating tableaux correspond to lattice paths. In the next
Section we provide some background on vacillating tableaux and the bijection.

\subsection{Vacillating tableaux}
A Young diagram (shape) is a collection of squares arranged in
left-justified
rows with weakly decreasing number of boxes in each row. A Young tableau
is a filling of the squares by numbers which is weakly decreasing
in each row and strictly decreasing in each column. A tableau is called
standard if each entry occurs exactly once. A tableau-sequence is a
sequence $\varnothing=\mu^{0}, \mu^{1},\ldots,\mu^{n}=
\varnothing$ of standard Young diagrams, such that for
$1\le i \le n$, $\mu^{i}$
is obtained from $\mu^{i-1}$ by either adding one square, removing
one square or doing nothing.
\begin{center}\label{F:pb}
\scalebox{0.45}[0.45]{\includegraphics*[20,710][600,810]{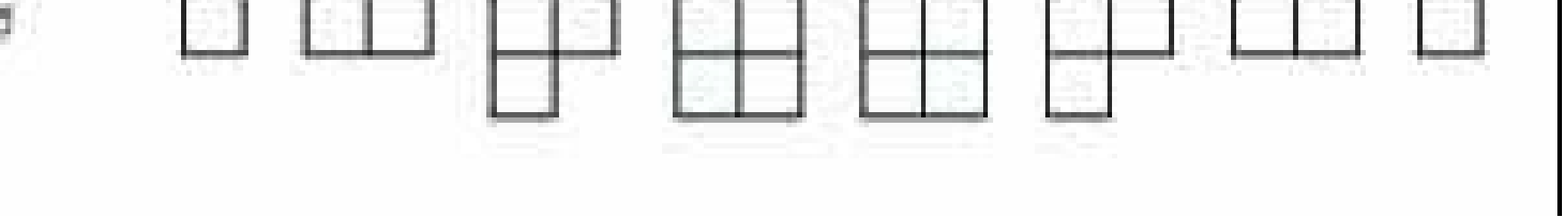}}
\end{center}
The RSK-algorithm is a process of row-inserting elements into a
tableau. Suppose we want to insert $k$ into a standard Young tableau
$\lambda$. Let $\lambda_{i,j}$ denote the element in the $i$-th row
and $j$-th column of the Young tableau. Let $i$ be the largest
integer such that $\lambda_{1,i-1}\le k$. (If $\lambda_{1,1}>k$,
then $i=1$.) If $\lambda_{1,i}$ does not exist, then simply add $k$
at the end of the first row. Otherwise, if $\lambda_{1,i}$ exists,
then replace $\lambda_{1,i}$ by $k$. Next insert $\lambda_{1,i}$
into the second row following the above procedure and continue until
an element is inserted at the end of a row. As a result we obtain a
new standard Young tableau with $k$ included. For instance inserting
the number sequence $5,2,4,1,6,3$ starting with an empty shape
yields the following sequence of standard Young tableaux:
\begin{center}\label{F:pb}
\scalebox{0.45}[0.45]{\includegraphics*[20,710][600,830]{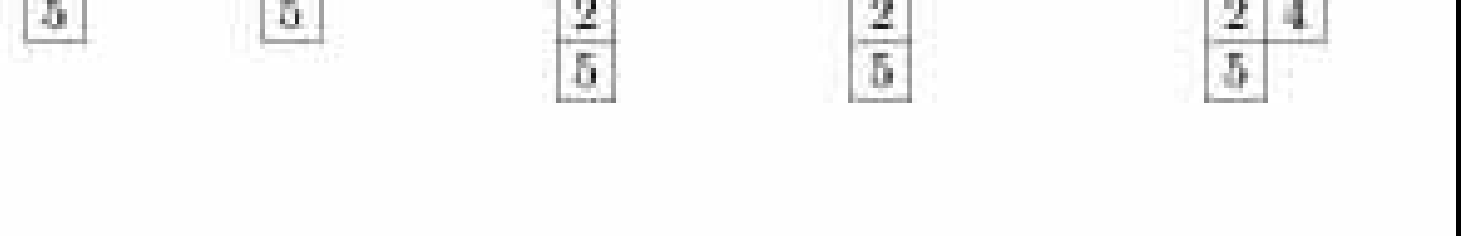}}
\end{center}
%%%%%%%%%%%%%%%%%%%%%%%%%%%%%%%%%%%%%%%%%%%%%%%%%%%%%%%%%%%%%%%%%%%%%%%%
%%%
A vacillating tableaux \cite{Reidys:07vac} $V_{\lambda}^{2n}$ of
shape $\lambda$ and length $2n$ is a sequence $(\lambda^{0},
\lambda^{1},\ldots, \lambda^{2n})$ of shapes such that {\sf (i)}
$\lambda^{0}=\varnothing$ and $\lambda^{2n}=\lambda,$ and {\sf (ii)}
$(\lambda^{2i-1},\lambda^{2i})$ is derived from
           $\lambda^{2i-2}$, for $1\le i\le n$ by either
$(\varnothing,\varnothing)$: doing nothing twice;
$(-\square,\varnothing)$: first removing a square then doing
nothing; $(\varnothing,+\square)$: first doing nothing then adding a
square; $(\pm \square,\pm \square)$: adding/removing a square at the
odd and even steps, respectively. Let $\mathcal{V}_\lambda^{2n}$
denote the set of vacillating tableaux.
%%%
%%%%%%%%%%%%%%%%%%%%%%%%%%%%%%%%%%%%%%%%%%%%%%%%%%%%%%%%%%%%%%%%%%%%%%%%
%%%
For instance, let us consider the following vacillating tableaux:
\begin{center}\label{F:pb}
\scalebox{0.45}[0.45]{\includegraphics*[20,710][600,840]{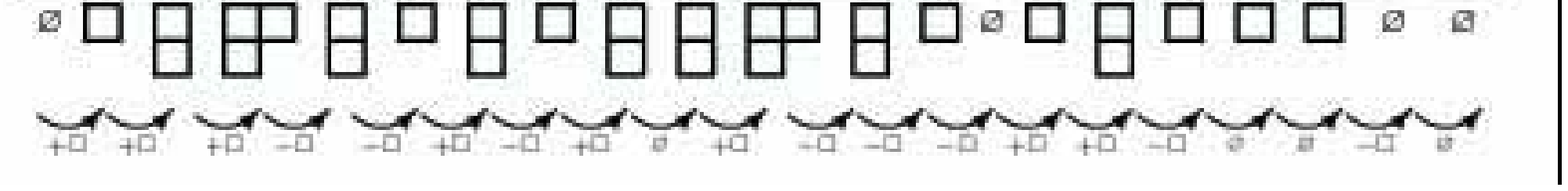}}
\end{center}
%%%%%%%%%%%%%%%%%%%%%%%%%%%%%%%%%%%%%%%%%%%%%%%%%%
\subsection{A bijection between vacillating tableaux and
tangled-diagrams}\label{S:bij}
When constructing the bijection
between vacillating tableaux and tangled-diagrams in
Theorem~\ref{T:bij} below, the notion of the inflation of a
tangled-diagram is important.
%%%
%%%%%%%%%%%%%%%%%%%%%%%%%%%%%%%%%%%%%%%%%%%%%%%%%%%%%%%%5
We are now able to discuss the bijection between vacillating
tableaux and tangled diagrams.
%%%
%%%%%%%%%%%%%%%%%%%%%%%%%%%%%%%%%%%%%%%%%%%%%%%%%%%%%%%%%%%%%%%%%%%%
%%%
\begin{theorem}\label{T:bij}\cite{Reidys:07vac}
There exists a bijection between the set of vacillating tableaux of shape
$\varnothing$ and length $2n$, $\mathcal{V}_\varnothing^{2n}$ and the
set of tangled-diagrams over $n$ vertices, $\mathcal{G}_n$
\begin{equation}
\beta\colon \mathcal{V}_{\varnothing}^{2n}  \longrightarrow
\mathcal{G}_n \ .
\end{equation}
Furthermore a tangled-diagram $G_n$ is $k$-noncrossing if and only if all
shapes $\lambda^i$ in its vacillating tableaux have less than $k$ rows.
That is $\phi\colon \mathcal{V}_\varnothing^{2n}\longrightarrow \mathcal{G}_n$
maps vacillating tableaux having less than $k$ rows into $k$-noncrossing
tangled-diagrams.
\end{theorem}
%%%
%%%%%%%%%%%%%%%%%%%%%%%%%%%%%%%%%%%%%%%%%%%%%%%%%%%%%%%%%%%%%%%%%%%
%%%
The proof of Theorem~\ref{T:bij} relies on the idea to resolve the
vertices of degree $2$ via an inflation, i.e.~vertex $i$ is resolved
by the pair $(i,i')$, where we utilize the linear order
$1<1'<2<2'<\dots <(n-1)<(n-1)'<n<n'$. The inflation transforms each
tangled-diagram into a partial matching with isolated points.
%%%%%%%%%%%%%%%%%%%%%%%%%%%%%%%%%%%%%%%%%%%%%%%%%%%%
%%%%%%%%%%%%%%%%%%%%%%%%%%%%%%%%%%%%%%%%%%%%%%%%%%%%%%%%%%%%%%%%
For instance,
\begin{figure}
\begin{center}
\scalebox{0.65}[0.65]{\includegraphics*[10,760][610,808]{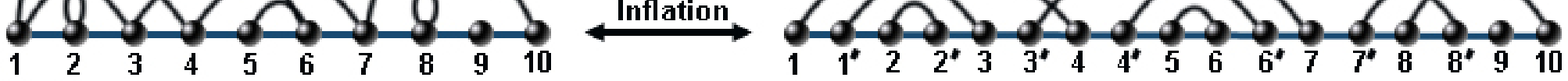}}
\end{center}
\caption{\small The inflation map: each vertex $i$ of degree $2$ is
replaced by a
pair of vertices, $(i,i')$, each incident to an respective arc.}
\end{figure}
%%%%%%
%%%%%%%%%%%%%%%%%%%%%%%%%%%%%%%%%%%%%%%%%%%%%%%%%%%%%%%%%%%%%%%%
%%%%%%
%%%
%%%%%%%%%%%%%%%%%%%%%%%%%%%%%%%%%%%%%%%%%%%%%%%%%%%%%%%%%%%%%%%%%%%%%%%%%%
\begin{figure}\label{F:bije}
\scalebox{0.6}[0.6]{\includegraphics*[20,580][610,850]{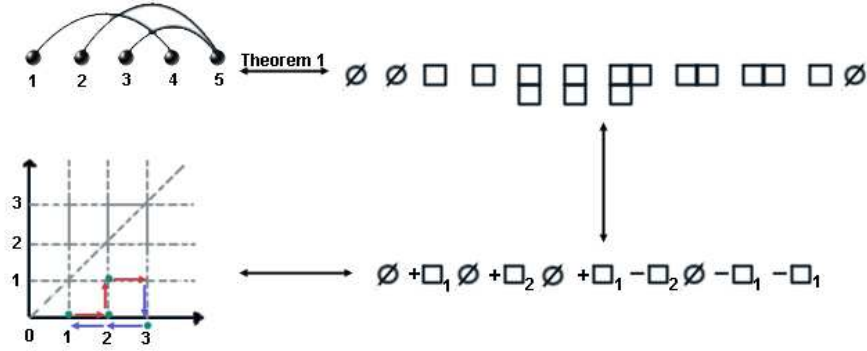}}
{\caption{\small From tangled-diagrams to lattice paths. First the
tangled-diagram (upper left) is resolved into its vacillating
tableaux (upper right). Reading the numbers of squares in the
corresponding rows (bottom right) induces the $2n$-step lattice path
(bottom right), which starts and ends in $(1,0)$. The path has
$\varnothing$ (green points), $+\square$ and $-\square$ (red and
purple points) induced by the pair steps $(\varnothing,+\square)$,
$(-\square,\varnothing)$ and $(-\square,-\square)$. Note that the
lattice path does not touch the ``wall'' $x=y$.}
 }
\end{figure}
%%%
%%%%%%%%%%%%%%%%%%%%%%%%%%%%%%%%%%%%%%%%%%%%%%%%%%%%%%%%%%%%%%%%%%%%%%%%
%%%
Restricting the steps for vacillating tableaux produces the bijections of
Chen {\it et.al} \cite{Chen}. Let $\mathcal{M}_k(n)$, $\mathcal{P}_k(n)$ and
$\mathcal{B}^{\dagger}_k(n)$ denote the set of $k$-noncrossing
matchings \cite{Stanley99}, partitions and braids without isolated
points over $[n]$, respectively. Theorem~~\ref{T:bij}
basically says the tableaux-sequences $\mathcal{M}_k(n)$,
$\mathcal{P}_k(n)$ and $\mathcal{B}^{\dagger}_k(n)$ are composed by
the elements in $S_{\mathcal{M}_k},S_{\mathcal{P}_k}$ and
$S_{\mathcal{B}^{\dagger}_k}$, respectively, where
\begin{eqnarray*}
S_{\mathcal{M}_k} & = & \{(-\square_h,\varnothing),(\varnothing,+\square_h)\}\\
S_{\mathcal{P}_k} & = & \{(-\square_h,\varnothing),(\varnothing,+\square_h),
(\varnothing,\varnothing),(-\square_h,+\square_l)\}\\
S_{\mathcal{B}^{\dagger}_k} & = &\{(-\square_h,\varnothing),
(\varnothing,+\square_h),(+\square_h,-\square_l)\} 1\le h,l\le k-1
\end{eqnarray*}
and $\pm \square_h$ denote the adding or subtracting of the rightmost square
``$\,\square_h\,$'' in the $h$th row in a given shape $\lambda$ and
let ``$\,\varnothing\,$'' denote doing nothing.
%%%
%%%%%%%%%%%%%%%%%%%%%%%%%%%%%%%%%%%%%%%%%%%%%%%%%%%%%%%%%%%%%%%%%%%%%%%%%%
%%%%%%%%%%%%%%%%%%%%%%%%%%%%%%%%%%%%%%%%%%%%%%%%%%%%%%%%%%%%%%%%%%%%%%%%
%%%
%%%
To get some intuition above the particular steps and
diagram-configurations let us show the key correspondences between
tableaux and diagram-motifs
\begin{center}\label{F:pb}
\scalebox{0.45}[0.45]{\includegraphics*[20,710][600,820]{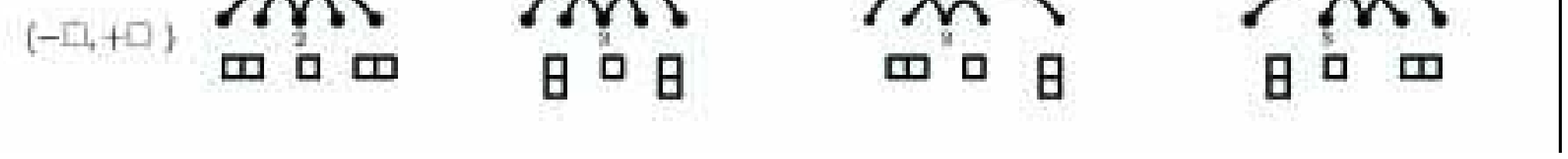}}
\end{center}

%%%
%%%%%%%%%%%%%%%%%%%%%%%%%%%%%%%%%%%%%%%%%%%%%%%%%%%%%%%%%%%%%%%%%%%%%%%%
%%%

\section{$k$-noncrossing tangled diagrams and $2$-regular, $k$-noncrossing
         partitions}

%%%
%%%%%%%%%%%%%%%%%%%%%%%%%%%%%%%%%%%%%%%%%%%%%%%%%%%%%%%%%%%%%%%%%%%%%%%%
%%%

In this section we prove two enumeration results. We give explicit
formulas for $k$-noncrossing tangled diagrams with a fixed number of
degree $2$ vertices and $2$-regular $k$-noncrossing partitions.
Since the latter formula is quite complicated we provide a simple
asymptotic expression in Section~\ref{S:Asym}.

Let $f_{k}(n)$ denote the number of perfect matching over $[n]$ and
$C_{m}$ be the Catalan number. Our first result reads

%%%
%%%%%%%%%%%%%%%%%%%%%%%%%%%%%%%%%%%%%%%%%%%%%%%%%%%%%%%%%%%%%%%%%%%%%%%%
%%%
\begin{theorem}
The number of the $k$-noncrossing tangled-diagrams over $[n]$ with
$\ell$ vertices of degree two, denoted by $d_{\ell,k}(n)$ is given
by
$$
d_{\ell,k}(n)=\sum_{i=0}^{n}{n \choose i}{n-i\choose
\ell}f_{k}(n-i+\ell)
$$
and in particular for $k=3$ we have
$$
d_{\ell,3}(n)=\sum_{i=0}^{n}{n \choose i}{n-i\choose
\ell}\left(C_{\frac{n-i+\ell}{2}}\, C_{\frac{n-i+\ell}{2}+2}-
C_{\frac{n-i+\ell}{2}+1}^2\right)
\ .\\
$$
\end{theorem}
%%%
%%%%%%%%%%%%%%%%%%%%%%%%%%%%%%%%%%%%%%%%%%%%%%%%%%%%%%%%%%%%%%%%%%%%%%%%
%%%
\begin{proof}
Let $\mathcal{D}_{i,\ell,k}$ denote the set of tangled-diagrams over
$[n]$ with $i$ isolated points and $\ell$ vertices of degree two and
$d_{i,\ell,k}=|\mathcal{D}_{i,\ell,k}|$. There are ${n \choose
i}{n-i \choose \ell}$ ways to choose the locations of the isolated
points and the vertices of degree two. Furthermore for an
arbitrary tangled-diagram over $V=[n]$ with $i$ isolated points
$V_1=\{v_1,v_2,\dots,v_i\}\subset V$ and $\ell$ vertices of degree
two $V_2=\{v_{i+1},\dots,v_{i+\ell}\}\subset V$, let
$\tilde{V}=V\setminus ( V_1\cup V_2) =\{v_{i+\ell+1},\dots,v_{n}\}$
be the set of vertices of degree one, via the inflation we will have
a perfect matching over $[|\{V_2\cup
V'_2\cup\tilde{V}\}|]=[2\ell+n-i-\ell]=[n-i+\ell]$, where
$V'_2=\{v'_{i+1},\dots,v'_{i+\ell}\}$. Since
$d_{\ell,k}=\sum_{i=0}^{n}d_{i,\ell,k}$, the theorem follows.
\end{proof}
%%%
%%%%%%%%%%%%%%%%%%%%%%%%%%%%%%%%%%%%%%%%%%%%%%%%%%%%%%%%%%%%%%%%%%%%%%%%
%%%
The first 10 number for $d_{i,\ell,3}$ for $\ell=1,2,3$ and
$n=1\dots 10$ are given by
\begin{center}
\begin{tabular}{|c|c|c|c|c|c|c|c|c|c|c|}
  \hline
  $\ell$,$n$ & 1 & 2 &3 & 4 & 5 & 6 & 7 & 8 & 9 & 10 \\
  \hline
  1 & 1 & 2 & 12 & 40 & 165 & 606 & 2380 & 9136 & 36099 & 142750 \\
  \hline
  2 & 0 & 3 & 9 & 102 & 450 & 2565 & 11823 & 57876 & 266220 & 1243170 \\
  \hline
  3 & 0 & 0 & 14 & 56 & 980 & 5320 & 38920 & 214144 & 1251852 & 6672120 \\
  \hline
\end{tabular}
\end{center}
%%%
%%%%%%%%%%%%%%%%%%%%%%%%%%%%%%%%%%%%%%%%%%%%%%%%%%%%%%%%%%%%%%%%%%%%%%%%
%%%
We proceed by enumerating $2$-regular $k$-noncrossing partitions. A
valid approach for this consists in building on the enumeration
results of \cite{MIRXIN} for $k$-noncrossing partitions using the
inclusion-exclusion principle. This strategy leads to functional
equations which prove that the asymptotic formulas of $2$-regular
$k$-noncrossing partitions and braids without isolated points
coincide. But braids can be enumerated via kernel methods
\cite{Knuth,Bander,Fayolle} directly, while $2$-regular
$k$-noncrossing partitions cannot. This suggests an alternative
ansatz \cite{Reidys:07dual}, by directly establishing a relation
between partitions and braids and consequently enumerating
partitions via braids. In Lemma~\ref{L:bijection} below we show this
correspondence. To this end we replace in a braid without isolated
points each loop by an isolated vertex and each pair of crossing
arcs at a degree $2$ vertex by noncrossing arcs, i.e.
\begin{center}\label{F:pb}
\scalebox{0.5}[0.5]{\includegraphics*[20,710][600,820]{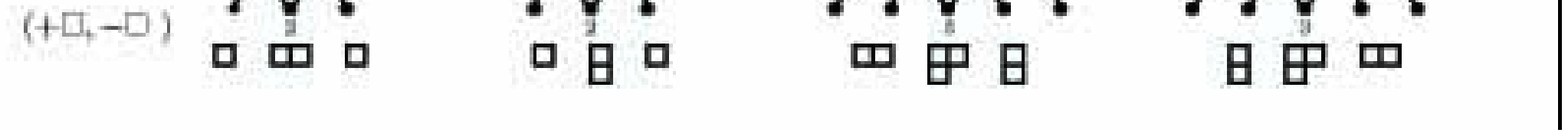}}
\end{center}
\begin{center}\label{F:pb}
\scalebox{0.7}[0.7]{\includegraphics*[20,710][600,820]{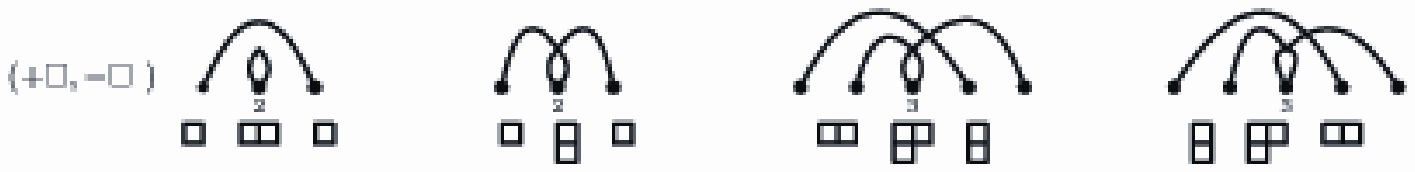}}
\end{center}
%%%
%%%%%%%%%%%%%%%%%%%%%%%%%%%%%%%%%%%%%%%%%%%%%%%%%%%%%%%%%%%%%%%%%%%%%%%%
%%%
Accordingly, we can identify braids without isolated points with a subset
of $3$-noncrossing partitions.

%%%
%%%%%%%%%%%%%%%%%%%%%%%%%%%%%%%%%%%%%%%%%%%%%%%%%%%%%%%%%%%%%%%%%%%%%%%%%%%%%
%%%
\begin{lemma}\label{L:bijection}
\cite{Reidys:07dual} Let $k\in\mathbb{N}$, $k\ge 3$. Then we have the
bijection
\begin{equation}\label{E:nobraid}
\vartheta\colon \mathcal{P}_{k,2}(n)\longrightarrow
\mathcal{B}_k^\dagger(n-1) \ ,
\end{equation}
where $\vartheta$ has the following property: for any
$\pi\in\mathcal{P}_k(n)$ holds: $(i,j)$ is an arc of $\pi$ if and
only if $(i,j-1)$ is an arc in $\vartheta(\pi)$.
\end{lemma}
\begin{proof}
By construction, $\vartheta$ maps tangled-diagrams over $[n]$ into
tangled diagrams over $[n-1]$. Since there exist no arcs of the form
$(i,i+1)$, $\vartheta(\pi)$ is, for any $\pi\in
\mathcal{P}_{k,2}(n)$ loop-free. By construction, $\vartheta$
preserves the orientation of arcs, whence $\vartheta(\pi)$
is a partition.\\
%%%%%%%%%%%%%%%%%%%%%%%%%%%%%%%%%%%%%%%%%%%%%%%%%%%%%%%%%%%%%%%%%%%%%%%%
{\it Claim.} $\vartheta\colon \mathcal{P}_{k,2}(n)\longrightarrow
\mathcal{B}_k^\dagger(n-1)$ is well-defined.\\
%%%%%%%%%%%%%%%%%%%%%%%%%%%%%%%%%%%%%%%%%%%%%%%%%%%%%%%%%%%%%%%%%%%%%%
We first prove that $\vartheta(\pi)$ is $k$-noncrossing. Suppose
there exist $k$ mutually crossing arcs, $(i_s,j_s)$, $s=1,\dots ,k$
in $\vartheta(\pi)$. Since $\vartheta(\pi)$ is a partition we have
$i_1<\dots <i_k<j_1<\dots <j_k$. Accordingly, we obtain for the
partition $\pi \in \mathcal{P}_{k,2}(n)$ the $k$ arcs $(i_s,j_s+1)$,
$s=1,\dots ,k$ where $i_1<\dots <i_k<j_1+1<\dots <j_k+1$, which is
impossible since $\pi$ is $k$-noncrossing. We next show that
$\vartheta(\pi)$ is a $k$-noncrossing braid. If $\vartheta(\pi)$ is
not a $k$-noncrossing braid, then according to eq.~(\ref{E:nobraid})
$\vartheta(\pi)$ contains $k$ arcs of the form $(i_1,j_1),\dots
(i_k,j_k)$ such that $i_1< \dots <i_k=j_1<\dots <j_k$ holds. Then
$\pi$ contains the arcs $(i_1,j_1+1)$, $(i_k,j_k+1)$ where
$i_1<\dots <i_k<j_1+1<\dots <j_k+1$, which is impossible since these
arcs are a set of $k$ mutually crossing arcs
and the claim follows.\\
%%%%%%%%%%%%%% inverse %%%%%%%%%%%%%%%%%%
{\it Claim.} $\vartheta$ is bijective.\\
Clearly $\vartheta$ is injective and it remains to prove
surjectivity. For any $k$-noncrossing braid $\delta$ there exists
some $2$-regular partition $\pi$ such that $\vartheta(\pi)=\delta$.
We have to show that $\pi$ is $k$-noncrossing. Let
$M'=\{(i_1,j_1),\dots ,(i_k,j_k)\}$ be a set of $k$ mutually
crossing arcs, i.e.~$i_1<\dots <i_k<j_1<\dots <j_k$. Then we have in
$\vartheta(\pi)$ the arcs $(i_s,j_s-1)$, $s=1,\dots ,k$ and
$i_1<\dots <i_k\le j_1-1<\dots <j_k-1$. If $M=\{(i_1,j_1-1),\dots
,(i_k,j_k-1)\}$ is $k$-noncrossing then we conclude $i_k=j_1-1$.
Therefore $M=\{(i_1,j_1-1),\dots ,(i_k,j_k-1)\}$, where $i_k=j_1-1$
which is, in view of eq.~(\ref{E:nobraid}) impossible in
$k$-noncrossing braids. By transposition we have thus proved that
any $\vartheta$-preimage is necessarily a $k$-noncrossing partition,
whence the claim and the proof of the lemma is complete.\\
\end{proof}
As an illustration of the bijection of Lemma~\ref{L:bijection} we display
%%%
%%%%%%%%%%%%%%%%%%%%%%%%%%%%%%%%%%%%%%%%%%%%%%%%%%%%%%%%%%%%%%%%%%%%%%%%
%%%
\begin{figure}
\begin{center}\label{F:pb}
\scalebox{0.55}[0.55]{\includegraphics*[10,700][530,750]{bp1.eps}}
\end{center}
\caption{\small The bijection $\vartheta\colon \mathcal{P}_{k,2}(n)
\longrightarrow \mathcal{B}_k^\dagger(n-1)$. Crossings are reduced
by contracting the arcs.}
\end{figure}
%%%
%%%%%%%%%%%%%%%%%%%%%%%%%%%%%%%%%%%%%%%%%%%%%%%%%%%%%%%%%%%%%%%%%%%%%%
%%%
Via Lemma~\ref{L:bijection} we have reduced the enumeration of
$2$-regular $k$-noncrossing partitions to that of braids without
isolated points. Let us discuss how the latter can be enumerated via
lattice paths. From Theorem \ref{T:bij}, (see Figure~\ref{F:bije})
we know that a $3$-noncrossing braid corresponds to a lattice paths
in the first
quadrant with the following properties:\\
{\sf(1)} the path starts and ends at $(1,0)$,\\
{\sf(2)} each step pair $(2i-1,2i)$, where $1\leq i\leq n$ is an
element of $$\{(0,+e_1), (0,+e_2),(-e_1,0), (-e_2,0),(+e_1,-e_1),
(+e_2,-e_2), (+e_1,-e_2),(+e_2,-e_1)\} \ .$$
{\sf(3)} the path never touches the wall $x=y$.\\
The key result facilitating the enumeration is the reflection
principle due D. Andr$\acute{e}$ in $1887$ \cite{andre} and
subsequently generalized by Gessel and Zeilberger \cite{Zeilberger}.
It is worth mentioning that this strategy is nonconstructive since
enumeration is obtained by counting {\it all} paths and having paths
touching the wall cancel each other.
%%%
%%%%%%%%%%%%%%%%%%%%%%%%%%%%%%%%%%%%%%%%%%%%%%%%%%%%%%%%%%%%%%%%%
%%%
\begin{theorem}{\bf (Reflection-Principle)}{\cite{Zeilberger}}\label{T:reflect}
Suppose $\mathcal{S}\in \{{\mathcal{M}_3}, { \mathcal{P}_3},
{\mathcal{B}_3}\}$ and let $\Omega_{\mathcal{S}}^{\text{\tiny \rm
(1,0)}}(2n)$ denote the number of $\mathcal{S}$-walks of length $2n$
from $(1,0)$ to $(1,0)$ that remain in the region $R=\{(x,y)\mid x>
y\ge 0,\ (x,y)\in \mathbb{Z}^2\}$. Let furthermore $f^{\text{\tiny
\rm $(x',y')$}}_{\text{\tiny \rm $(x,y)$}} (2n)$ be the number of
$\mathcal{S}$-walks from $(x,y)$ to $(x',y')$ of length $2n$ that
remain in the first quadrant. Then we have
\begin{equation}\label{E:refl}
\Omega_{\mathcal{S}}^{\text{\tiny \rm (1,0)}}(2n) = f^{\text{\tiny
\rm $(1,0)$}}_{\text{\tiny \rm $(1,0)$}}(2n)- f^{\text{\tiny \rm
$(0,1)$}}_{\text{\tiny \rm $(1,0)$}}(2n) \ .
\end{equation}
\end{theorem}
%%%
%%%%%%%%%%%%%%%%%%%%%%%%%%%%%%%%%%%%%%%%%%%%%%%%%%%%%%%%%%%%%%%%%%%%%%%%
%%%
\begin{proof}
Suppose $\gamma$ is a $\mathcal{S}$-walk starting and ending at
$(1,0)$ which remains in the first quadrant and that touches the
diagonal $x=y$. Let $(a,a)$ be the first point where $\gamma$
touches the diagonal $y=x$. Reflect all steps of $\gamma$ after
$\gamma$ touched the diagonal in $(a,a)$ and denote the resulting
walk by $\gamma'$. Then $\gamma'$ is a $\mathcal{S}$-walk starting
from $(1,0)$ and ending at $(0,1)$. This procedure yields a unique
pair $(\gamma,\gamma')$ for each $\mathcal{S}$-walk $\gamma$
starting and ending at $(1,0)$ which remains in the first quadrant
and that touches the diagonal $x=y$. According to eq.~(\ref{E:refl})
these pairs cancel themselves and only the paths that never touch
the diagonal remain, whence the theorem.
\end{proof}
%%%
%%%%%%%%%%%%%%%%%%%%%%%%%%%%%%%%%%%%%%%%%%%%%%%%%%%%%%%%%%%%%%%%%%5
%%%
\begin{figure}
\begin{center}
\scalebox{0.55}[0.55]{\includegraphics*[20,550][600,808]{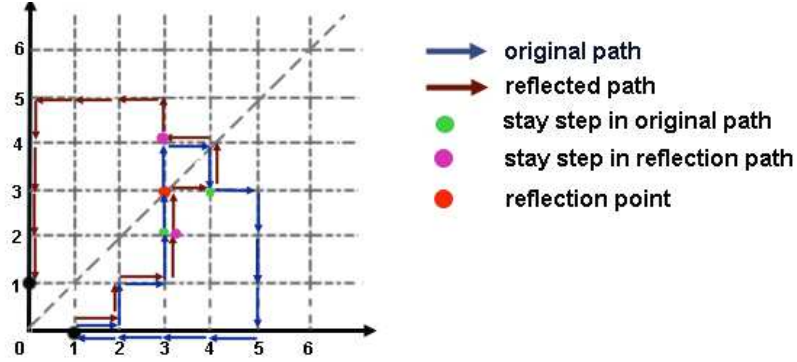}}
\end{center}
\caption{The reflection principle. The original lattice path(blue)
starting and ending at $(1,0)$ touches the wall $x=y$ at $(3,3)$ for
the first time. The corresponding reflected path(red) starts at
$(1,0)$ and ends at $(0,1)$ obtained by reflecting all steps after
$(3,3)$ w.r.t.~the wall $x=y$.}
\label{F:pb}
\end{figure}
%%%
%%%%%%%%%%%%%%%%%%%%%%%%%%%%%%%%%%%%%%%%%%%%%%%%%%%%%%%%%%%%%%%%%%
%%%
%%%
%%%%%%%%%%%%%%%%%%%%%%%%%%%%%%%%%%%%%%%%%%%%%%%%%%%%%%%%%%%%%%%%%%%%%%%%
%%%
Using the reflection principle we can enumerate braids via the
kernel method \cite{Knuth,Bander,Fayolle}. In fact these computation
have been obtained by \cite{MIRXIN} who enumerated enhanced
partitions. Our second result reads
%%%
%%%%%%%%%%%%%%%%%%%%%%%%%%%%%%%%%%%%%%%%%%%%%%%%%%%%%%%%%%%%%%%%%%%%%%%5
%%%
\begin{theorem}\label{T:braid}
The number of $2$-regular, $3$-noncrossing partitions is given by
\begin{eqnarray*}{\label{E:braidnum1}}
p_{3,2}(n+1) &=&
\sum_{s\in\mathbb{Z}}\left[\beta_{n}(1,0,s)-\beta_{n}(1,-1,s)-
\beta_{n}(1,-4,s)+\beta_{n}(1,-3,s)\right.\\
& & -\beta_{n}(3,4,s)+\beta_{n}(3,3,s)+\beta_n(3,0,s)-\beta_{n}(3,1,s)\\
& & \left.
+\beta_{n}(2,5,s)-\beta_{n}(2,4,s)-\beta_{n}(2,1,s)+\beta_{n}(2,2,s))\right]
\ ,
\end{eqnarray*}
where $\beta_{n}(k,m,s)=\frac{k}{n+1}{n+1\choose s}{n+1\choose k+s}
{n+1\choose s+m}$. Furthermore $p_{3,2}(n)$ satisfies the recursion
\begin{equation}\label{E:recursion-b}
\alpha_1(n)\, p_{3,2}(n+1)+ \alpha_2(n)\, p_{3,2}(n+2)+
\alpha_3(n)\, p_{3,2}(n+3)- \alpha_4(n)\, p_{3,2}(n+4)=0 \ ,
\end{equation}
where
\begin{eqnarray*}
\alpha_1(n) & = & 8(n+2)(n+3)(n+1) \\
\alpha_2(n) & = & 3(n+2)(5n^2+47n+104) \\
\alpha_3(n) & = & 3(n+4)(2n+11)(n+7) \\
\alpha_4(n) & = & (n+9)(n+8)(n+7)  \ .
\end{eqnarray*}
\end{theorem}
%%
%%%%%%%%%%%%%%%%%%%%%%%%%%%%%%%%%%%%%%%%%%%%%%%%%%%%%%%%%%%%%%%%%%%%%%%%%%%%%
%%%
For instance, the first 12 numbers of $2$-regular, $3$-noncrossing
partitions are given by
$$
\begin{tabular}{|c|c|c|c|c|c|c|c|c|c|c|c|c|c|}
  \hline
  $n$ & 1 & 2 & 3 & 4 & 5 & 6 & 7 & 8 & 9 & 10 & 11 &12\\
  \hline
  $p_{3,2}(n)$ & 1 &1& 2 & 5 & 15 & 51 & 191 & 772 & 3320 & 15032 & 71084
&348889\\
  \hline
\end{tabular}
$$
We will show in the next section that the formulas
given in Theorem~\ref{T:braid} have simple asymptotic
formulas.

\section{Asymptotic analysis}\label{S:Asym}

%%%%
%%%%%%%%%%%%%%%%%%%%%%%%%%%%%%%%%%%%%%%%%%%%%%%%%%%%%%%%%%%%%%%%%%%%%%%%
%%%

In this
section we employ the particularly elegant theory of singular
difference equations due to Birkhoff and Trjitzinsky \cite{T:wimp}.
The theory of Birkhoff-Trjitzinsky establishes form, existence and
properties of such fundamental sets in general, and will be
discussed in the Appendix. For our purposes it suffices to identify
the unique, monotonously increasing formal series solution (FSS).
%%%%%%%%%%%%%%%%%%%%%%%%%%%%%%%%%%%%%%%%%%%%%%%%%%%%%%%%%%%%%%%%%%%%%%%%
%%%
\begin{theorem}\label{T:2}
There exists some real constants $K>0$ and $c_{1},c_{2},c_{3}$ such that
\begin{equation}
p_{3,2}(n+1)\sim K \ 8^{n}n^{-7}(1+c_{1}/n+c_{2}/n^2+c_3/n^3)
\end{equation}
holds. Explicitly, we have $K=6686.408973$, $c_1=-28,\
c_2=455.77778$ and $c_3=-5651.160494$.
\end{theorem}
%%%%
%%%%%%%%%%%%%%%%%%%%%%%%%%%%%%%%%%%%%%%%%%%%%%%%%%%%%%%%%%%%%%%%%%%%%%%%
%%%
\begin{proof}
{\emph{Claim.}} There exists some $K>0$ and $c_{1},c_{2},c_{3}\dots$
such that
\begin{equation}
p_{3,2}(n+1)\sim K \ 8^{n}n^{-7}(1+c_{1}/n+c_{2}/n^2+c_3/n^3\cdots).
\end{equation}
 Theorem~\ref{T:Birk} guarantees the existence of $3$
linearly independent formal series solutions (FSS) for
eq.~(\ref{E:recursion-b}). We proceed by constructing these using
the following ansatz for $p_{3,2}(n)$:
\begin{equation}
p_{3,2}(n+1) = E(n)K(n) \quad E(n)=e^{\mu_0n\ln
n+\mu_{1}n}n^{\theta}
\end{equation}
where
\begin{equation}
K(n)=\exp\{\alpha_{1}n^{\beta+\alpha_2n^{\beta-1/\rho+\cdots}}\},\ \
\ \ \alpha_1\neq 0,\ \beta=j/\rho,\ 0\leq j<\rho.
\end{equation}
We immediately derive setting $\lambda =e^{\mu_{0}+\mu_{1}}$
\begin{align*}
\frac{p_{3,2}(n+k+1)}{p_{3,2}(n+1)}&=n^{\mu_{0}k}\lambda^{k}
\{1+\frac{k\theta+k^2\mu_{0}/2}{n}+\cdots\} \\
  & \quad \ \exp\{\alpha_1\beta kn^{\beta-1}+\alpha_2(\beta-\frac{1}{\rho})
kn^{\beta-1/\rho-1+\cdots}\}.
\end{align*}
We arrive at
\begin{align*}
0=1+&\frac{15}{8}\{1+\frac{\theta+\mu_0/2+\frac{27}{5}}{n}+
\cdots\}\xi\{1+(\alpha_1\beta
n^{\beta-1}+\alpha_2(\beta-1/\rho)n^{\beta-1/\rho-1}+\cdots)+\cdots\}\\
+&\frac{3}{4}\{1+\frac{2\theta+2\mu_0+\frac{21}{2}}{n}+\cdots\}
\xi^2\{1+(2\alpha_1\beta
n^{\beta-1}+2\alpha_2(\beta-1/\rho)n^{\beta-1/\rho-1}+\cdots)+\cdots\}\\
-&\frac{1}{8}\{1+\frac{3\theta+9\mu_0/2+18}{n}+\cdots\}
\xi^3\{1+(3\alpha_1\beta
n^{\beta-1}+3\alpha_2(\beta-1/\rho)n^{\beta-1/\rho-1}+\cdots)+\cdots\}.
\end{align*}
First we consider the maximum power of $n$, which is zero. In view
of $1=\frac{1}{8}n^{3\mu_0}\lambda^3$ we obtain $\mu_0=0$. This
implies $\rho=1$ since $\rho\geq 1$ and $\rho$ should be the
smallest integer s.t. $\rho\mu_0\in \mathbb{N}$. Equating the
constant terms again, we obtain that $\lambda$ is indeed a root of
the cubic polynomial $P(X)$
$$
P(X) = 1+\frac{15}{8}X+\frac{3}{4}X^2-\frac{1}{8}X^3.
$$
Therefore we have $\lambda=8$ or $-1$. Notice that $0\leq \beta <1$
implies $\beta=0$. Otherwise, equating the coefficient of
$n^{\beta-1}$ implies $\alpha_1=0$, which is impossible. It remains
to compute $\theta$. For this purpose we equate the coefficient of
$n^{-1}$, i.e.~ $ 8\frac{15}{8}(\theta+\frac{27}{5})+8^2\frac{3}{4}
(\frac{21}{2}+2\theta)-8^3\frac{1}{8}(18+3\theta)=0 $ from which we
conclude $\theta=-7$. Since $p_{3,2}(n)$ is monotone increasing
$p_{3,2}(n)$ coincides with the only monotonously increasing FSS,
given by
\begin{equation}
p_{3,2}(n+1)\sim K \cdot 8^{n}\cdot
n^{-7}(1+c_{1}/n+c_{2}/n^2+c_3/n^3\cdots)
\end{equation}
for some $K>0$ and constants $c_1,c_2,c_3$ and the proof of the
claim is complete. We compute $c_1=-28$, $c_2=455.778$ and
$c_3=-5651.160494$ by equating the coefficients of $n^{-2},\ n^{-3}$
and $n^{-4}$, ($2268+81c_1=0,\ 1683c_1+162c_2-26712=0$ and
$-32547c_1+729c_2+129654+243c_3=0$) and finally get $K=6686.408973$
numerically to complete the proof of the theorem.
\end{proof}
%%%%%%%%%%%%%%%%%%%%%%%%%%%%%%%%%%%%%%%%%%%%%%%%%%%%%%%%%%%%%%%%%%%%%%%%
%%%
%%%
%%%%%%%%%%%%%%%%%%%%%%%%%%%%%%%%%%%%%%%%%%%%%%%%%%%%%%%%%%%%%%%%%%%%%%%%
%%%
\section{Appendix}
%%%%%%%%%%%%%%%%%%%%%%%%%%%%%%%%%%%%%%%
\subsection{The Birkhoff-Trjitzinsky theory} Any difference
equation with rational coefficients \cite{wimp} can be written as
\begin{equation}\label{E:1}
\sum_{h=0}^m C_m(n) \, y(n+h) = 0 \ \quad C_0(n)=1,\quad C_m(n)\neq 0,
\quad n=0,1,2,\dots
\end{equation}
where the coefficients possess representations as generalized Poincar\'e
series
\begin{equation}\label{E:3}
C_h(n)\sim n^{\frac{K_h}{\omega}}\left[ c_{0,h}+c_{1,h}n^{-\frac{1}{\omega}}
+ c_{2,h}n^{-\frac{2}{\omega}}+ \dots \right], \quad h=1,2,\dots \ .
\end{equation}
Here $K$ is an integer, $\omega$ is an integer $\ge 1$ independent of $h$
and $c_{0,h}\neq 0$ unless $C_h(n)=0$.
We shall assume that $\omega$ is minimal.
A set of functions $z^{(j)}(n)$ is called linearly independent if the
determinant
\begin{equation}
\forall \; n \in \mathbb{N}\cup \{0\}; \quad
\det{(z^{(j+1)}(n+i))}_{0\le i,j\le h-1} \neq 0.
\end{equation}
The classical theory of difference equations asserts that eq.~(\ref{E:1})
possesses a set of linearly independent solutions constituting a basis of
the solution space. Such a set is called a fundamental set. The
Birkhoff-Trjitzinsky theory proves
that there exists a fundamental set in which all elements have an asymptotic
expansion consisting of an exponential leading term multiplied by a linear
combination of descending series of the form eq.~(\ref{E:3}).
To provide the notion of {\it formal series solution} and {\it Birkhoff series}
we set
\begin{eqnarray*}
Q(\rho,n) & = &
\mu_0 n \ln(n) +\sum_{j=1}^\rho \mu_j n^{\frac{\rho+1-j}{\rho}}, \qquad
s(\rho,n)  =  n^\theta\, \sum_{j=0}^t (\ln(n))^j n^{\frac{r_{t-j}}{\rho}}
q_j(\rho,n), \\
q_j(\rho,n) & = & \sum_{s=0}^\infty b_{sj} n^{-\frac{s}{\rho}}
\end{eqnarray*}
where $\rho,r_j,\mu_0\rho$ are integers, $\rho\ge 1$, $\mu_j,\theta,b_{sj}\in
\mathbb{C}$, $b_{0,j}\neq 0$, unless  $b_{sj}=0$ for $s=0,1,2,\dots$, $r_0=0$,
$-\pi\le \text{\rm Im}(\mu_1)<\pi$. Then we call
$$
y(\rho,n)= e^{Q(r,n)}s(\rho,n)
$$
a formal series solution (FSS) of eq.~(\ref{E:1}) if and only if substituted
in eq.~(\ref{E:1}) after dividing by $e^{Q(\rho,n)}$ and corresponding
algebraic transformations, the coefficients of
$$
n^{\theta+\frac{r}{\rho}+\frac{s}{\omega}}\ln(n)^j, \quad r,s=0,1,
\dots,t \quad r,s=0,\pm 1,\pm 2, \dots,
$$
are equal to zero. For given sequence $(f(n))_{n\ge 0}$ we furthermore call
\begin{equation}
f(n) \sim e^{Q(\rho,n)}s(\rho,n)
\end{equation}
the Birkhoff series for $f(n)$ if and only if for every $k\ge 1$ there exist
bounded functions $A_{kj}(n)$, $j=0,1,\dots,t$, such that
\begin{equation}
e^{-Q(\rho,n)} n^{-\theta} f(n) =
\sum_{j=0}^t\ln(n)^jn^{\frac{r_{t-j}}{\rho}}
\sum_{s=0}^{k-1} b_{sj} n^{\frac{s}{\rho}}+n^{-\frac{k}{\rho}}
\sum_{j=0}^t \ln(n)^jn^{\frac{r_{t-j}}{\rho}} A_{kj}(n) \ .
\end{equation}
Following \cite{wimp} we define
\begin{equation}
w_k= \det{(e^{Q_{j+1}(\rho,n+i)}s_{j+1}(\rho,n+i))}_{0\le i,j\le k-1} \ .
\end{equation}

The main result of the Birkhoff-Trjitzinsky theory can now be stated as
follows
%%%
%%%%%%%%%%%%%%%%%%%%%%%%%%%%%%%%%%%%%%%%%%%%%%%%%%%%%%%%%%%%%%%%%%%%%%%%%
%%%
\begin{theorem}\cite{Birkhoff,T:wimp}\label{T:Birk}
There exist exactly $m$ FSS of eq.~(\ref{E:1}) of type
$e^{Q(\rho,n)}s(\rho,n)$ where $\rho=\nu\omega$ for some integer
$\nu\ge 1$ and each FSS represents asymptotically some solution of
the equation. The above FSS are, up to multiplicative constants,
unique and the $m$ solutions so represented constitute a fundamental
set for the equation.
\end{theorem}
%%%
%%%%%%%%%%%%%%%%%%%%%%%%%%%%%%%%%%%%%%%%%%%%%%%%%%%%%%%%%%%%%%%%%%%%%%%%%
%%%

%%%
%%%%%%%%%%%%%%%%%%%%%%%%%%%%%%%%%%%%%%%%%%%%%%%%%%%%%%%%%%%%%%%%%%%%%%%%%
%%%

%%%
%%%%%%%%%%%%%%%%%%%%%%%%%%%%%%%%%%%%%%%%%%%%%%%%%%%%%%%%%%%%%%%%%%%%%%%%%%
%%%
{\bf Acknowledgments.}
%%%
%%%%%%%%%%%%%%%%%%%%%%%%%%%%%%%%%%%%%%%%%%%%%%%%%%%%%%%%%%%%%%%%%%%%%%%%%%
%%%
We are grateful to Emma Y. Jin for helpful discussions.
This work was supported by the 973 Project, the PCSIRT Project of the
Ministry of Education, the Ministry of Science and Technology, and
the National Science Foundation of China.

%%%
%%%%%%%%%%%%%%%%%%%%%%%%%%%%%%%%%%%%%%%%%%%%%%%%%%%%%%%%%%%%%%%%%%%%%%%%
%%%

\bibliographystyle{amsplain}

\begin{thebibliography}{10}

\bibitem{Science:05}\textit{Mapping RNA form and function}. Science, \textbf{2}, 2005.
\bibitem{andre} D. Andr$\acute{e}$ \textit{Solution directed du proble$\grave{e}$me, r$\acute{e}$solu par M. Bertrand}. C. R. Acad. Sci. Paris \textbf{105}(1887)
436--437.
\bibitem{NMR} Alexander S. Brodsky, Heidi A. Erlacher and James R.
Williamson, \textit{NMR evidence for a base triple in the HIV-2 TAR}
$\textbf{C-G}\cdot \textbf{C}^{+}$ \textit{mutant-argininamide
complex.} Nucleic Acids Research, \textbf{26} (1998), No. 8,
1991--1995.
\bibitem{MIRXIN} Mireille Bousquet-M\'{e}lou and Guoce Xin,
\textit{On partitions avoiding 3-crossings.} S\'{e}minaire
Lotharingien de Combinatoire, \textbf{54} (2006), Article B54c.

\bibitem{Birkhoff} George D. Birkhoff, \textit{Formal theory of irregular difference equations} Acta Math. \textbf{54}
(1930), 205-246.

\bibitem{T:wimp} George D. Birkhoff and W. J. Trjitzinsky, \textit{Analytic theory of singular difference equations} Acta Math.,  \textbf{60}
(1932), 1--89.

\bibitem{Chen} William Y. C. Chen, Eva Y.P. Deng, Rosena R.X. Du,
Richard P. Stanley and Catherine H. Yan, \textit{Crossings and
Nestings of Matchings and Partitions.} Trans. Amer. Math. Soc.
\textbf{359}  (2007), No. 4, 1555--1575.
\bibitem{Reidys:07vac} William Y. C. Chen, Jing Qin and Christian M.
Reidys, \textit{Crossings and Nestings of tangled-diagrams.} Submitted.

\bibitem{Chastain} Michael Chastain and Ignacio Tinoco, Jr.,
\textit{A Base-triple Structural Domain in RNA} Biochemistry,
\textbf{31} (1992), 12733-12741.
\bibitem{Batey}Robert T. Batey, Robert P. Rambo, and Jennifer A.
Doudna, \textit{Tertiary Motifs in RNA Structure and Folding} Angew.
Chem. Int. Ed., \textbf{38} (1999), 2326-2343.

\bibitem{Bander} C. Banderier, M. Bousquet-M\'{e}lou, A. Denise, P.
Flajolet, D. Gardy, and D. Gouyou-Beuchamps, \textit{Generating
functions of generating trees}, Discrete mathematics, \textbf{246}
(2002), no. 1-3, 29-55.

\bibitem{berele} A. Berele, \textit{A Schensted-type correspondence for the
  symplectic group}, J. Combinatorial Theory (A), \textbf{43}
  (1986), 320--328.
\bibitem{Fayolle}G. Fayolle, R. Iasnogorodske, and V. Malyshev.,
\textit{Random walks in the quarter-plane: Algebraic methods,
boundary value problems and applications,} volume 40 of
\textit{Applications of Mathematics.} Springer-Verlag, Berlin,1999.

\bibitem{Zeilberger}. M. Gessel and D. Zeilberger,  Random walk in a Weyl chamber,
\textit{Proc. Amer. Math. Soc.} \textbf{115} (1992), 27--31.

\bibitem{GB89}D. Gouyou-Beauschamps, \textit{Standard Young tableaux
of height  4 and 5}, Europ.\ J. Combin., \textbf{10} (1989), 69--82.

\bibitem{GM93}D. Grabiner and P. Magyar, \textit{Random walks in Weyl chambers
  and the decomposition of tensor powers}, J. Alg.\ Combinatorics,
  \textbf{2} (1993),  239--260.

\bibitem{Greene74}C. Greene, \textit{An extension of Schensted's theorem},
Adv.\  Math., \textbf{14} (1974),  254--265.

\bibitem{Stadler:99}C. Haslinger  and P.F. Stadler,\textit{ RNA
Structures with Pseudo-Knots. Bull.Math.Biol.}, \textbf{61}
 (1999)  437--467.

\bibitem{Reidys:07pseu}E.Y. Jin, J. Qin, and C.M. Reidys, \textit{
Combinatorics of RNA structures with pseudoknots}, Bull.Math.Biol.,
2007. in press.
\bibitem{Reidys:07dual}E.Y. Jin, J. Qin, and C.M. Reidys, \textit{On $k$-noncrossing partitions}, submited.
\bibitem{Reidys:07asym}E.Y. Jin and C.M. Reidys, \textit{Asymptotic enumeration of RNA structures with Pseudoknots},
Bull.Math.Biol., 2007. in press.

\bibitem{Knuth}D.E. Knuth, \textit{The art of computer programming},
vol. 1: Fundamental Algorithms, Addison-Wesley, 1973, Third edition,
1997.
\bibitem{Konings}D.A.M. Konings and R.R. Gutell, \textit{A comparison of thermodynamic foldings with comparatively derived structures of 16s and 16s-like rRNAs}, RNA,  \textbf{1} (1995), 559-574.

\bibitem{Loria}A. Loria and T. Pan, \textit{Domain structure of the ribozyme from eubacterial ribonuclease p}, RNA,  \textbf{2} (1996), 551-563.

\bibitem{Mohanty79} S. G. Mohanty, \textit{Lattice Path Counting and
    Applications}, Academic Press, New York, 1979.
\bibitem{Westhof}E. Westhof and L. Jaeger \textit{RNA pseudoknots}, Current Opinion Struct. Biol.,  \textbf{2} (1992), 327-333.
\bibitem{Wilf} M. Petov\v{s}ek, H.S. Wilf and D. Zeiberger.
\textit{A=B}. A K Peter Ltd., Wellesey, MA, 1996.

\bibitem{Sch61} C. E. Schensted, \textit{Longest increasing and decreasing
subsequences}, Canad.\ J. Math.,  \textbf{13} (1961), 179--191.
\bibitem{Shen} L. X. Shen and Zhuoping Cai and Ignacio Tinoco.
Jr,\textit{RNA structure at high resolution}, FASEB J., \textbf{9}
(1995),
 1023--1033.
\bibitem{Sivakova} Sona Sivakova and Stuart J. Rowan, \textit{Nucleobases
as supramolecular motifs} Chem. Soc. Rev.,  \textbf{34} (2005),
9-21.

\bibitem{ec1} R. Stanley, \textit{Enumerative Combinatorics}, vol.~1,
  Wadsworth and Brooks/Cole, Pacific Grove, CA, 1986;
  second printing, Cambridge University Press, Cambridge, 1996.

\bibitem{Stanley99} R. Stanley, \textit{Enumerative Combinatorics},
  vol.~2, Cambridge University Press, Cambridge, 1999.

\bibitem{wimp} Jet Wimp and Doron Zeilberger, \textit{Resurrecting the asymptotics of linear recurrences}, Journal of Mathmatical analysis and
applications,  \textbf{3} (1985), 162--176.

\bibitem{Zuker:81} M. Zuker and P. Stiegler, \textit{Optimal computer folding of large RNA sequence using thermodynamics and auxiliary
informations}, Nucl. Acid Res., \textbf{9} (1981), 133--148.

\bibitem{Zuker:84} M. Zuker and D. Sankoff, \textit{RNA Secondary Structure and their prediction}, Bull. Math. Biol., \textbf{46} (1984), 591-621.
\end{thebibliography}

\end{document}